\documentclass[12pt]{article}

\newcommand{\real}{\mathbb{R}}

\renewcommand{\natural}{\mathbb{N}}

\newcommand{\card}{{\rm card}}

\newcommand{\loc}{{\rm loc}}
\newcommand{\e}{{\rm e}}                %% \e = 2.771828...
\renewcommand{\d}{\,{\rm d}}            %% integration measure
\newcommand{\D}{{\rm d}}                %% differential operator
\def\1{{\bf 1}}                         %% \I = \sqrt{-1}
\newcommand{\erfc}{{\rm erfc}}
\renewcommand{\epsilon}{\varepsilon}

\newcommand{\CC}{{\cal C}}
\newcommand{\DD}{{\cal D}}
\newcommand{\EE}{{\cal E}}

\newcommand{\OO}{{\cal O}}

\newtheorem{theorem}{Theorem}[section]
\newtheorem{lemma}[theorem]{Lemma}

\newtheorem{proposition}[theorem]{Proposition}
\newtheorem{corollary}[theorem]{Corollary}
\newtheorem{remark}[theorem]{Remark}

\newcommand{\reff}[1]{(\ref{#1})}

\newcommand{\proof}{{\noindent \bf Proof:\ }}

\def\build#1_#2^#3{\mathrel{
  \mathop{\kern 0pt#1}\limits_{#2}^{#3}}}
\def\QED{\mbox{}\hfill$\Box$}

\newdimen\texpscorrection
\newdimen\figcenter
\def\figurewithtex #1 #2 #3 #4 #5\cr{\null
  {\goodbreak\figcenter=\hsize\relax
  \advance\figcenter by -#4truecm
  \divide\figcenter by 2
  \begin{figure}[hbt]
  \vskip #3truecm\noindent\hskip\figcenter
  \includegraphics{#1}{\hskip\texpscorrection\input #2 }
  \vskip 0.8truecm{\baselineskip=0.8\baselineskip
  \noindent \vbox{\noindent {\footnotesize #5}}\par}
  \end{figure}}}
\def\point#1 #2 #3 {\rlap{\kern #1 truecm
\raise #2 truecm \hbox{#3}}}

\usepackage{amsfonts}
\begin{document}

\title{A variational proof of global stability for bistable 
travelling waves}

\author{{\bf Thierry Gallay} \\ 
Institut Fourier\\
Universit\'e de Grenoble I\\
38402 Saint-Martin-d'H\`eres\\ 
France\\
\and
{\bf Emmanuel Risler}\\
Institut Camille Jordan\\
INSA de Lyon\\
69621 Villeurbanne\\
France}
%\date{December 22, 2006}

\maketitle
\begin{abstract}
We give a variational proof of global stability for bistable
travelling waves of scalar reaction-diffusion equations on the real
line. In particular, we recover some of the classical results by
P.~Fife and J.B.~McLeod (1977) without any use of the maximum
principle. The method that is illustrated here in the simplest
possible setting has been successfully applied to more general
parabolic or hyperbolic gradient-like systems.
\end{abstract}

%%%%%%%%%%%%%%%%%%%%%%%%%%%%%%%%%%%%%%%%%%%%%%%%%%%%%%%%%%%%%%%%%%%%%%%%%%%

\section{Introduction} \label{intro}

The purpose of this work is to revisit the stability theory for
travelling waves of reaction-diffusion systems on the real line. We
are mainly interested in {\em global} stability results which assert
that, for a wide class of initial data with a specified behavior at
infinity, the solutions approach for large times a travelling wave
with nonzero velocity. In the case of scalar reaction-diffusion
equations, such properties have been established by Kolmogorov,
Petrovski \& Piskunov \cite{kpp}, by Kanel \cite{kanel1, kanel2}, and
by Fife \& McLeod \cite{fifemcleod1, fifemcleod2} under various
assumptions on the nonlinearity.  The proofs of all these results use
a priori estimates and comparison theorems based on the parabolic
maximum principle. Therefore they cannot be extended to general
reaction-diffusion systems nor to scalar equations of a different
type, such as damped hyperbolic equations or higher-order parabolic
equations, for which no maximum principle is available.  However,
these methods have been successfully applied to monotone
reaction-diffusion systems \cite{RTV,volpert}, as well as to scalar
equations on infinite cylinders \cite{roquejoffre1, roquejoffre2}.

Recently, a different approach to the global stability of bistable
travelling waves has been developped by the second author
\cite{rislerglob}. The new method is of variational nature and is
therefore restricted to systems which admit a gradient structure, but
it does not make any use of the maximum principle and is therefore
potentially applicable to a wide class of problems. The goal of 
this paper is to explain how this method works in the simplest 
possible case, namely the scalar parabolic equation
\begin{equation}\label{equ}
  u_t \,=\, u_{xx} - F'(u)~,
\end{equation}
where $u = u(x,t) \in \real$, $x \in \real$, and $t \ge 0$. We shall
thus recover the main result of Fife \& McLeod \cite{fifemcleod1}
under slightly different assumptions on the nonlinearity $F$, with a
completely different proof. The present article can also serve as an
introduction to the more elaborate work \cite{rislerglob}, where the
method is developped in its full generality and applied to the
important case of gradient reaction-diffusion systems of the form $u_t
= u_{xx} - \nabla V(u)$, with $u \in \real^n$ and $V : \real^n \to
\real$. A further application of our techniques is given in 
\cite{gallayjoly}, where the global stability of travelling 
waves is established for the damped hyperbolic equation 
$\alpha u_{tt} + u_t = u_{xx} - F'(u)$, with $\alpha > 0$.

We thus consider the scalar parabolic equation \reff{equ}, 
which models the propagation of fronts in chemical reactions
\cite{billinghamneedham}, in combustion theory \cite{kanel1,kanel2}, 
and in population dynamics \cite{aronsonweinberger,fisher}. We
suppose that the ``potential'' $F : \real \to \real$ is a smooth, 
coercive function with a unique global minimum and at least one
additional local minimum. More precisely, we assume that $F \in 
\CC^2(\real)$ satisfies
\begin{equation}\label{coercive}
  \liminf_{|u|\to\infty} uF'(u) \,>\, 0~.
\end{equation}
In particular, $F(u) \to +\infty$ as $|u| \to \infty$. We also assume 
that $F$ reaches its global minimum at $u = 1$: 
\begin{equation}\label{min1}
  F(1) \,=\, -A < 0~, \quad F'(1) \,=\, 0~, \quad  F''(1) \,>\, 0~,  
\end{equation}
and has in addition a local minimum at $u = 0$:
\begin{equation}\label{min0}
  F(0) \,=\, F'(0) \,=\, 0~, \quad  F''(0) \,=\, \beta > 0~.  
\end{equation}
Finally, we suppose that all the other critical values of $F$ 
are positive, namely
\begin{equation}\label{other}
  \Bigl\{u \in \real \,\Big|\, F'(u) = 0\,,~F(u) \le 0\Bigr\} 
  \,=\, \{0\,;1\}~.
\end{equation}
A typical potential satisfying the above requirements is represented
in Fig.~1. 

\figurewithtex 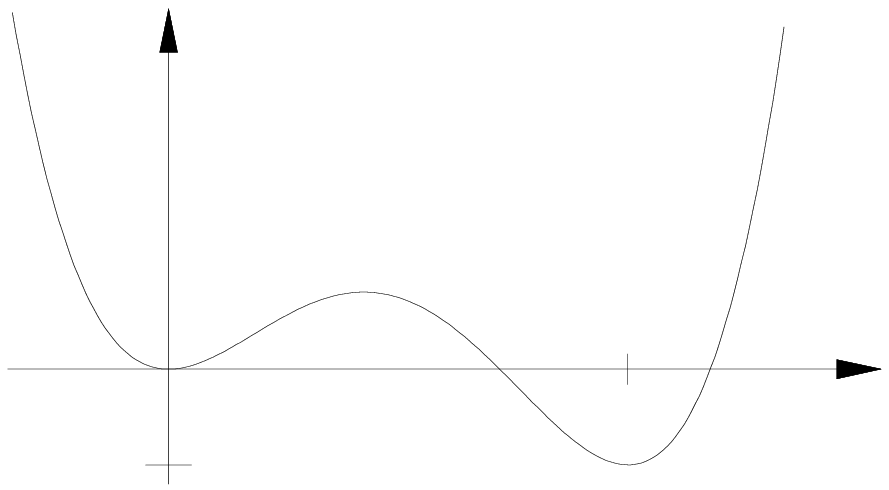 Fig1.tex 5.000 9.000
{\bf Fig.~1:} The simplest example of a nonlinearity $F$ satisfying 
assumptions \reff{coercive}--\reff{other}.\cr

Under assumptions \reff{min1}-\reff{other}, it is well-known that
Eq.\reff{equ} has a family of travelling waves of the form $u(x,t) =
h(x -c_*t)$ connecting the stable equilibria $u = 1$ and $u = 0$. More
precisely, there exists a unique speed $c_* > 0$ such that the
boundary value problem
\begin{equation}\label{hdef}
\left\{\begin{array}{l}
  h''(y) + c_* h'(y) - F'(h(y)) \,=\, 0~, \quad y \in \real~,\\
  h(-\infty) = 1~, \quad h(+\infty) = 0~,
\end{array}\right.
\end{equation}
has a solution $h : \real \to (0,1)$, in which case the profile 
$h$ itself is unique up to a translation. Moreover $h \in
\CC^3(\real)$, $h'(y) < 0$ for all $y \in \real$, and $h(y)$ 
converges exponentially to its limits as $y \to \pm\infty$. 

This family of travelling waves plays a major role in the dynamics 
of Eq.\reff{equ}, as is shown by the following global convergence 
result:

\begin{theorem}\label{thm1}
Let $F \in \CC^2(\real)$ satisfy assumptions \reff{coercive}--\reff{other}. 
Then there exist $\delta > 0$ and $\nu > 0$ such that, for all initial 
data $u_0 \in \CC^0(\real)$ with
\begin{equation}\label{lims}
 \limsup_{x\to -\infty} |u_0(x)-1| \,\le\, \delta~, \quad
 \limsup_{x\to +\infty} |u_0(x)| \,\le\, \delta~,
\end{equation}
Eq.\reff{equ} has a unique global bounded solution satisfying $u(x,0) 
= u_0(x)$ for all $x \in \real$. In addition, there exists 
$x_0 \in \real$ such that
\begin{equation}\label{conc}
  \sup_{x\in \real} \Big|u(x,t) - h(x - c_*t -x_0)\Big| \,=\, 
  \OO(\e^{-\nu t})~, \quad \hbox{as }t \to +\infty~.
\end{equation}
\end{theorem}

Theorem~\ref{thm1} was first proved by Fife \& McLeod
\cite{fifemcleod1,fifemcleod2} under the additional assumption that $0
\le u_0(x) \le 1$ for all $x \in \real$. In that case $u(x,t) \in
[0,1]$ for all $x \in \real$ and all $t \ge 0$ by the maximum
principle, so that the coercivity assumption \reff{coercive} is not
needed.  As is mentioned in \cite{fife}, the results of
\cite{fifemcleod1} can be extended to arbitrary initial data
satisfying \reff{lims} provided that $u F'(u) > 0$ for all $u \notin
[0,1]$, a condition that is more restrictive than \reff{coercive} in
the sense that $F$ is not allowed to have critical points outside the
interval $[0,1]$. The simplest case considered in \cite{fifemcleod1}
is when $F$ has exactly one critical point in the open interval
$(0,1)$, a situation in which condition \reff{other} is clearly met.
However, Fife \& McLeod also study the case where $F$ has three
critical points in the open interval, including a local minimum at $u
= u_* \in (0,1)$. In this situation there exists a travelling wave
solution of \reff{equ} with speed $c_1 > 0$ connecting $u = 1$ to $u =
u_*$, and also a travelling wave with speed $c_2 \in \real$ connecting
$u = u_*$ to $u = 0$. If $c_1 > c_2$, which is always the case if
\reff{other} holds, there exists $c_* \in (0,1)$ such that \reff{hdef}
has a solution $h : \real \to (0,1)$, and the conclusion of
Theorem~\ref{thm1} is still valid. If $c_1 < c_2$, there exists no
travelling wave connecting $u = 1$ to $u = 0$, and the solution of
\reff{equ} with initial data satisfying \reff{lims} converges as $t
\to \infty$ to a superposition of two travelling waves
\cite{fifemcleod1}.

Theorem~\ref{thm1} is a particular case of the general results
obtained in \cite{rislerglob}, see Theorem~4 in Section~9.6 of that
reference. Therefore, there is no need to give here a complete proof.
Instead we shall prove the convergence result \reff{conc} under the
additional assumption that the initial data $u_0(x)$ decay rapidly to
zero as $x \to +\infty$. It is intuitively clear that the precise
behavior of $u_0(x)$ near $x = +\infty$ should not play an important
role, because the equilibrium $u = 0$ ahead of the front is stable
(this is in sharp contrast with the case of a monostable front invading
an unstable equilibrium, where the behavior ahead of the front is of
crucial importance). However, this restriction allows to shortcut many
technicalities and to give a much simpler proof in which the essence
of the argument can be easily understood.

Our approach is based on the fact that Eq.\reff{equ} possesses (at 
least formally) a {\em gradient structure}, not only in the laboratory
frame but also in any frame moving to the right with a positive 
velocity. To see this, we introduce the following notation. If
$u(x,t)$ is a solution of \reff{equ}, we define for any $c > 0$
\begin{equation}\label{vdef}
  v(y,t) \,=\, u(y+ct,t)~, \quad \hbox{or equivalently}\quad  
  u(x,t) \,=\, v(x-ct,t)~. 
\end{equation}
Setting $y = x-ct$ we see that the new function $v(y,t)$ satisfies
\begin{equation}\label{eqv}
  v_t \,=\, v_{yy} + c v_y - F'(v)~.
\end{equation}
We now introduce the {\em energy functional}
\begin{equation}\label{EEcdef}
  \EE_c[v] \,=\, \int_\real e^{cy}\Bigl(\frac12\,v_y^2 + F(v)\Bigr)
  \d y~,
\end{equation}
and the corresponding {\em energy dissipation functional}
\begin{equation}\label{DDcdef}
  \DD_c[v] \,=\, \int_\real e^{cy} \Bigl(v_{yy} + c v_y -F'(v)
  \Bigr)^2\d y~.
\end{equation}
We also denote by $H^1_c(\real)$ the Banach space
\begin{equation}\label{H1cdef}
  H^1_c(\real) \,=\, \Bigl\{v \in L^\infty(\real)\,\Big|\, e^{cy/2}v 
  \in H^1(\real)\Bigr\}~,
\end{equation}
equipped with the norm $\|v\|_{H^1_c} = \|v\|_{L^\infty} + 
\|e^{cy/2}v\|_{H^1}$. Note that any $v \in H^1_c(\real)$ decays to 
zero faster than $e^{-cy/2}$ as $y \to +\infty$. Since $F(v) \sim
\beta v^2/2$ as $v \to 0$ by \reff{min0}, it follows that $\EE_c[v]
< \infty$ for all $v \in H^1_c(\real)$. Conversely, any $v \in 
L^\infty(\real)$ such that $v(y) \to 0$ as $y \to +\infty$ belongs
to $H^1_c(\real)$ as soon as $\EE_c[v] < \infty$. 

If $v(y,t)$ is a solution of \reff{eqv} with initial data $v_0 \in 
H^1_c(\real)$, then $v(\cdot,t) \in H^1_c(\real)$ for all $t \ge 0$
and a direct calculation shows that
\begin{equation}\label{ED}
  \frac{\D }{\D t}\,\EE_c[v(\cdot,t)] \,=\, -\DD_c[v(\cdot,t)] 
  \,\le\, 0~, \quad t > 0~. 
\end{equation}
In other words, the energy $\EE_c$ is a {\em Lyapunov function} of
system \reff{eqv} in $H^1_c(\real)$. This observation is of course not
new: in their original proof, Fife \& McLeod \cite{fifemcleod1}
already used a suitable truncation of the functional $\EE_c$ for the
particular value $c = c_*$ to show that the solution $v(y,t)$ of
\reff{eqv} approaches a travelling wave for a sequence of times.
However, the fact that Eq.\reff{equ} has a whole family of
(nonequivalent) Lyapunov functions has not been fully exploited until
recently. The only reference we know where the implications of this
rich Lyapunov structure are really discussed is a recent paper by
Muratov \cite{muratov}, which contains a lot of interesting
observations and a few general results concerning a wider class of
systems than Eq.\reff{equ}, but fails to prove the convergence to
travelling waves. The goal of the present article is to show that, in
the simple case of Eq.\reff{equ}, the gradient structure {\em alone}
is sufficient to establish convergence, at least if we restrict
ourselves to solutions which decay to zero rapidly enough as $x \to
+\infty$ so that the energy functionals are properly defined.

The main difficulty of this purely variational approach is that we do
not have good a priori estimates on the solution $v(y,t) = u(y+ct,t)$
in a moving frame with speed $c > 0$. First of all, it is not clear a
priori that the energy $\EE_c[v(\cdot,t)]$ is bounded from below
(this will not be the case typically if $c$ is too small), and without
this information it is difficult to really exploit the dissipation
relation \reff{ED}. Next, if we have a lower bound on
$\EE_c[v(\cdot,t)]$, we can deduce from \reff{ED} that the solution
$v(y,t)$ converges uniformly on compact sets, at least for a sequence
of times, towards a stationary solution of \reff{eqv}, but we cannot
exclude a priori that this limit is just the trivial equilibrium $v
\equiv 0$ (this will be the case typically if $c$ is too large). To
overcome these difficulties, the main idea is to track the position of
the front interface in the following way. We fix positive constants
$\beta_1, \beta_2$ such that $\beta_1 < F''(0) < \beta_2$, and we
choose $\epsilon > 0$ small enough so that
\begin{equation}\label{epsdef}
  \beta_1 \,\le\, F''(u) \,\le\, \beta_2~,
  \quad \hbox{for all } u \in [-2\epsilon,2\epsilon]~.
\end{equation}
Given a continuous solution of \reff{equ} satisfying the boundary 
conditions
\begin{equation}\label{limits}
  \lim_{x \to -\infty}u(x,t) \,=\, 1~, \quad 
  \lim_{x \to +\infty}u(x,t) \,=\, 0~, \quad t \ge 0~,
\end{equation}
we define the {\em invasion point} $\bar x(t)$ as the first point 
starting from the right where the solution $u(x,t)$ leaves an 
$\epsilon$-neighborhood of the equilibrium $u = 0$:
\begin{equation}\label{barxdef}
  \bar x(t) \,=\, \max\Bigl\{x \in \real~\Big|~ |u(x,t)| \ge 
  \epsilon \Bigr\}~.
\end{equation}
In view of \reff{limits}, it is clear that $-\infty < \bar x(t) <
\infty$ for all $t \ge 0$, and that $|u(\bar x(t),t)| = \epsilon$.  
A quantity similar to $\bar x(t)$ was also introduced in \cite{muratov},
where it is called the ``leading edge''.

The strategy of the proof is to show that the solution $u(x,t)$ 
converges uniformly on compact sets around the invasion point 
$\bar x(t)$ towards a suitable translate of the travelling wave
\reff{hdef}. Using only the gradient structure, we can prove the
following result: 

\begin{proposition}\label{thm2}
Let $F \in \CC^2(\real)$ satisfy assumptions \reff{coercive}--\reff{other}.
If $u_0 \in H^1_c(\real)$ for some sufficiently large $c > 0$ and
$u_0 - 1 \in H^1(\real_-)$, then the solution $u(x,t)$ of 
Eq.\reff{equ} with initial data $u_0$ satisfies, for all $L > 0$, 
\begin{equation}\label{locconv}
  \sup_{z \in [-L,+\infty)}|u(\bar x(t)+z,t) - h_\epsilon(z)| 
  \,\build\hbox to 8mm{\rightarrowfill}_{t \to \infty}^{}\, 0~,
\end{equation}
where $\bar x(t)$ is the invasion point \reff{barxdef} and
$h_\epsilon$ is the travelling wave \reff{hdef} normalized so that
$h_\epsilon(0) = \epsilon$. Moreover the map $t \mapsto \bar x(t)$ 
is $\CC^1$ for $t$ sufficiently large and $\bar x'(t) \to c_*$ 
as $t \to \infty$.
\end{proposition}

As is explained above, the assumption $u_0 \in H^1_c(\real)$ is 
needed in order to use the energy functional $\EE_c$ without
truncating the unbounded exponential factor $e^{cy}$. The proof
will show that is it sufficient to take here $c > \sqrt{2A}/\epsilon$,
where $A$ is defined in \reff{min1} and $\epsilon$ in \reff{epsdef}. 
On the other hand, the assumption $u_0 - 1 \in H^1(\real_-)$ 
is just a convenient way to guarantee that the first condition in 
\reff{limits} is satisfied, but with minor modifications we can 
treat the more general case where $|u_0(x)-1|$ is assumed to be small
for large $x < 0$, as in \reff{lims}. 

The local convergence established in Proposition~\ref{thm2} is the key
step in proof of Theorem~\ref{thm1}. Once \reff{locconv} is known, it
remains to show that the solution $u(x,t)$ converges uniformly to $1$
in the region far behind the invasion point $\bar x(t)$. Such a 
``repair'' is certainly expected because $u = 1$ is the point where
the potential $F$ reaches its global minimum. A convenient way to
prove this is to use a truncated version of the functional
\begin{equation}\label{EEdef}
  \EE[u] \,=\, \int_\real \Bigl(\frac12\,u_x^2 + \overline{F}(u)\Bigr)
  \d x~,
\end{equation}
where $\overline{F}(u) = F(u) - F(1) \ge 0$. In this way, we can show
that the solution $u(x,t)$ approaches uniformly on $\real$ a
travelling wave (at least for a sequence of times), and using in
addition the local stability results established in \cite{sattinger}
we obtain \reff{conc}. We thus have:

\begin{corollary}\label{thm3}
Under the assumptions of Proposition~\ref{thm2}, there exist 
$x_0 \in \real$ and $\nu > 0$ such that \reff{conc} holds. 
\end{corollary}

We conclude this introduction with a few comments on the scope of 
our method. First, it is clear that the assumptions 
\reff{coercive}--\reff{other} are not the weakest ones under which
Proposition~\ref{thm2} holds. A careful examination of the proof
reveals that the only hypotheses that we really use are:

\begin{itemize}

\item[{\bf H1:}] {\em For all bounded initial data $u_0$, 
Eq.\reff{equ} has a (unique) global bounded solution.} This is 
certainly true if \reff{coercive} holds, but it is sufficient to 
assume, for instance, that $F(u) \to +\infty$ as $|u| \to \infty$, 
or that $uF'(u) > 0$ whenever $|u|$ is sufficiently large. 

\item[{\bf H2:}] {\em $F(0) = F'(0) = 0$, and there exists 
$\epsilon > 0$ such that $F''(u) \ge 0$ for all $u \in [-\epsilon,
\epsilon]$.} This is automatically true if \reff{min0} holds, 
but $u = 0$ need not be a strict local minimum of $F$. In 
particular Proposition~\ref{thm2} holds for the nonlinearities 
of combustion type considered in \cite{kanel1,kanel2}. 

\item[{\bf H3:}] {\em There exists a unique $c > 0$ such that 
the differential equation $v_{yy} + cv_y - F'(v) = 0$ has a 
bounded solution satisfying $|v(0)| = \epsilon$, $|v(y)| \le 
\epsilon$ for all $y \ge 0$, and $v(y) \to 0$ as $y \to +\infty$;
furthermore, this solution is unique.} Under assumptions 
\reff{min1}--\reff{other}, we have $c = c_*$ and $v = h_\epsilon$. 
In general, we can assume without loss of generality that 
that $v$ is positive and converges to $1$ as $y \to -\infty$, 
so that $F(1) < 0$ and $F'(1) = 0$. It also follows that 
$F(u) \ge 0$ for all $u \le 0$ and that $F$ has no critical point 
$u_* < 1$ with $F(u_*) < 0$.

\end{itemize}
On the other hand, to prove that the solution of \reff{equ}
given by Proposition~\ref{thm2} converges uniformly on $\real$ to 
a travelling wave we need the additional assumption:

\begin{itemize}

\item[{\bf H4:}] {\em There exists $\epsilon' > 0$ such that the 
only bounded solution of the differential equation $u_{xx} - F'(u) = 0$
with $|u(0)-1| \le \epsilon'$ is $u \equiv 1$.} This requires that
$F$ attains its global minimum at $u = 1$, and nowhere else. 

\end{itemize}
Finally, if we want the convergence to be exponential in time as 
in \reff{conc}, we need to assume that $F''(1) > 0$.

\medskip
Another comment concerns the variational structure of Eq.\reff{equ}.
Due to the exponential weight $e^{cy}$, it is clear that the 
energy functional $\EE_c$ is not translation invariant. In fact, 
for any $v \in H^1_c(\real)$ and any $\ell \in \real$, we have the
relation $\EE_c[v(\cdot - \ell)] = e^{c\ell}\EE_c[v]$. This implies 
that the infimum of $\EE_c[v]$ is either $0$ or $-\infty$. Under 
our assumptions on $F$, the transition between both regimes occurs 
precisely at the critical speed $c_*$ for which travelling waves
exist:
$$
  \inf_{v \in H^1_c(\real)} \EE_c[v] \,=\, \left\{
  \begin{array}{ccc} 0 & \hbox{if} & c \ge c_*~, \\
  -\infty & \hbox{if} & c < c_*~.
  \end{array}\right.
$$
Indeed, as was observed by Muratov \cite{muratov}, for any 
$c < c_* + \sqrt{c_*^2 + 4F''(0)}$ we have the identity
$$
  c \,\EE_c[h] \,=\, (c-c_*) \int_\real e^{cy} h'(y)^2 \d y~,
$$
where $h$ is the solution of \reff{hdef}. This shows in particular
that $\EE_c[h] < 0$ when $c < c_*$, hence $\inf \EE_c = -\infty$
in that case. The fact that $\EE_c \ge 0$ when $c \ge c_*$ is not
obvious a priori, and will be established in the course of the 
proof of Proposition~\ref{thm2}, see Corollary~\ref{conv3}. Note
also that $\EE_{c_*}[h] = 0$, so that $\inf \EE_{c_*} = \min 
\EE_{c_*} = 0$. 

The rest of the paper is organized as follows. In
Section~\ref{prelim}, we establish the basic inequalities relating 
the energy $\EE_c$, the dissipation $\DD_c$, and the invasion point. 
Using these relations, we prove in Section~\ref{invspeed} that
the average speed of the invasion point $\bar x(t)$ has a 
limit $c_\infty > 0$ as $t \to \infty$. The core of the paper 
is Section~\ref{localconv}, where we show that $c_\infty = c_*$ 
and prove Proposition~\ref{thm2}. The proof of Corollary~\ref{thm3}
is then performed in the final Section~\ref{repair}. 

\noindent{\bf Acknowledgements.} The authors are indebted to 
S. Heinze, R. Joly, and C.B. Muratov for fruitful discussions.

\section{Preliminary estimates} \label{prelim}

As the potential $F$ is smooth and coercive, it is well-known that
the Cauchy problem for the semilinear equation \reff{equ} is 
globally well-posed in the space of bounded functions, see 
e.g. \cite{henry}. Due to parabolic regularization, the solutions
are smooth for $t > 0$ and satisfy \reff{equ} in the classical 
sense. Under assumption \reff{coercive}, one can also show that
our system has a {\em bounded absorbing set} in the following 
sense:

\begin{lemma}\label{bas} There exists a constant $B > 0$ depending 
only on $F$ such that, for all initial data $u_0 \in L^\infty(\real)$, 
the (unique) solution $u(x,t)$ of \reff{equ} satisfies, 
for all sufficiently large $t \ge 0$, 
\begin{equation}\label{Bdef}
  \sup_{x\in\real}\Bigl(|u(x,t)|+|u_x(x,t)|+|u_{xx}(x,t)|\Bigr) 
  \,\le\, B~.
\end{equation}
Moreover, $u(\cdot,t)$ is bounded in $H^s_\loc(\real)$ for 
some $s > 5/2$ and all $t \ge 1$.
\end{lemma}

The uniform bound on $|u(x,t)|$ follows easily from the maximum 
principle, but it can also be established using localized energy 
estimates, see \cite[Section~9.1]{rislerglob}. The bounds  
on the derivatives are then obtained in a standard way using 
parabolic regularization.

From now on, we suppose that $u_0 \in H^1_{c_0}(\real)$ for some 
$c_0 > 0$ (which will be specified later) and that $u_0 - 1 \in 
H^1(\real_-)$. Then the solution of \reff{equ} with initial data
$u_0$ satisfies $u(\cdot,t) \in H^1_{c_0}(\real)$ and $u(\cdot,t) 
- 1 \in H^1(\real_-)$ for all $t \ge 0$, because $u = 0$ and $u = 1$
are (stable) equilibria of \reff{equ}. In particular, the boundary 
conditions \reff{limits} hold for all times, so that one can 
define the invasion point $\bar x(t)$ by \reff{barxdef}. Also, 
since we are interested in the long-time behavior of $u(x,t)$, 
we can assume without loss of generality that estimate \reff{Bdef}
is valid for all $t \ge 0$. 

As is explained in the introduction, we shall use the energy
functionals $\EE_c$ (for various values of $c > 0$) to prove that the
solution $u(x,t)$ converges to a travelling wave $h$ locally around
the invasion point $\bar x(t)$. A technical problem we shall encounter
is that the invasion point, as defined in \reff{barxdef}, need not be
a continuous function of time and can therefore jump back and forth in
an uncontrolled way. It is possible to avoid this difficulty using a
more clever definition than \reff{barxdef}, see \cite{rislerglob},
but we follow here another approach and just introduce a second
invasion point defined by

\begin{equation}\label{barXdef}
  \bar X(t) \,=\, \max\Bigl\{x \in \real~\Big|~ |u(x,t)| \ge 
  2\epsilon \Bigr\}~.
\end{equation}
Clearly, $-\infty < \bar X(t) < \bar x(t) < +\infty$ for all $t \ge 0$. 
The important point is that an information on $\bar x$ at a given time 
provides an upper bound on $\bar X$ at later times:

\begin{lemma}\label{Xxcontrol}
There exists $T_0 > 0$ and $C_0 > 0$ such that, for all $t_0 \ge 0$,
one has
\begin{equation}\label{Xxineg}
  \bar X(t) \,\le\, \bar x(t_0)+C_0 \quad \hbox{for all } t \in 
  [t_0,t_0+T_0]~.
\end{equation}
\end{lemma}

\proof Fix $t_0 \ge 0$. The solution of \reff{equ} satisfies
$$
  u(t) \,=\, S(t-t_0)u(t_0) - \int_{t_0}^t S(t-s)F'(u(s))\d s
  \,\equiv\, u_1(t) + u_2(t)~, \quad t \ge t_0~, 
$$
where $S(t) = e^{t\partial_x^2}$ is the heat semigroup. Take 
$K > 0$ such that $|F'(u)| \le K$ whenever $|u| \le B$, where
$B$ is as in \reff{Bdef}. Then $\|u_2(t)\|_{L^\infty} \le K(t-t_0)$.
On the other hand, by definition of $\bar x$, we have $|u(x,t_0)| 
\le \epsilon$ if $x \ge \bar x(t_0)$ and $|u(x,t_0)| \le B$ if 
$x \le \bar x(t_0)$. Using the explicit form of the heat kernel, we 
deduce that 
$$
  |u_1(x,t)| \,\le\, \frac{1}{\sqrt{4\pi (t-t_0)}}\int_{\real}
  e^{-\frac{(x-y)^2}{4(t-t_0)}} |u(y,t_0)|\d y \,\le\, \epsilon \,+\,
  \frac{B}2 \,\erfc\Bigl(\frac{x-\bar x(t_0)}{\sqrt{4(t-t_0)}}\Bigr)~,
$$
where $\erfc(x) = (2/\sqrt{\pi})\int_x^\infty e^{-z^2}\d z$.
We first choose $T_0 > 0$ such that $K T_0 < \epsilon/2$, 
and then $C_0 > 0$ such that $B\,\erfc(C_0/\sqrt{4T_0}) < \epsilon$. 
Then, for all $t \in [t_0,t_0+T_0]$ and all $x \ge \bar x(t_0) 
+ C_0$ we have $|u(x,t)| < 2\epsilon$, which implies that 
$\bar X(t) \le \bar x(t_0) + C_0$. \QED

\medskip
We now derive the basic estimates on the energy \reff{EEcdef} and
the energy dissipation \reff{DDcdef} which will be used throughout
the proof. Given $c \in (0,c_0)$, we define $v(y,t) = u(y+ct,t)$ 
as in \reff{vdef}, and we set
\begin{equation}\label{EDdef}
  E_c(t) \,=\, \EE_c[v(\cdot,t)]~, \quad 
  D_c(t) \,=\, \DD_c[v(\cdot,t)]~, \quad t \ge 0~.
\end{equation}
Of course, $v(y,t)$ depends also on the speed parameter $c$, but
to simplify the notations this dependence will not be indicated
explicitly. We also denote by $\bar y_c(t)$, $\bar Y_c(t)$ the
invasion points in the moving frame: 
\begin{equation}\label{barydef}
  \bar y_c(t) \,=\, \bar x(t) - ct~, \qquad 
  \bar Y_c(t) \,=\, \bar X(t) - ct~~. 
\end{equation}
By construction, $|v(\bar y_c(t),t)| = \epsilon$, $|v(y,t)| \le 
\epsilon$ for $y \ge \bar y_c(t)$ and $|v(y,t)| \le 2\epsilon$ 
for $y \ge \bar Y_c(t)$. Remark that, by \reff{epsdef}, the 
following inequalities hold whenever $|v| \le 2\epsilon$:
\begin{equation}\label{farfield}
  \frac{\beta_1}2\,v^2 \,\le\, F(v) \,\le\, \frac{\beta_2}2\, v^2~, 
  \quad \beta_1 v^2 \,\le\, v F'(v) \,\le\, \beta_2 v^2~,\quad
  \beta_1 \,\le\, F''(v) \,\le\, \beta_2~.
\end{equation}

\medskip\noindent{\bf Lower bound on $E_c\,$:} Using \reff{EEcdef}, 
\reff{farfield}, and the fact that $F(u) \ge -A$ for all $u \in
\real$, we find
\begin{eqnarray}\nonumber
  E_c(t) &=& \int_{-\infty}^{\bar y_c(t)}
  e^{cy}\Bigl(\frac12\,v_y^2(y,t) + F(v(y,t))\Bigr) \d y + 
  \int_{\bar y_c(t)}^\infty e^{cy}\Bigl(\frac12\,v_y^2(y,t) + 
  F(v(y,t))\Bigr) \d y \\ \label{lowprelim}
  &\ge& \int_{-\infty}^{\bar y_c(t)} e^{cy} (-A)\d y + 
  \int_{\bar y_c(t)}^\infty e^{cy}\Bigl(\frac12\,v_y^2(y,t) + 
  \frac{\beta_1}2\,v^2(y,t)\Bigr) \d y~.
\end{eqnarray}
To estimate the last integral in \reff{lowprelim} we recall that 
$v(\bar y_c(t),t)^2 = \epsilon^2$, so that
$$
  e^{c\bar y_c(t)}\epsilon^2 \,=\, -\int_{\bar y_c(t)}^\infty 
  \partial_y \Bigl(e^{cy}v^2(y,t)\Bigr)\d y \,=\, 
  -\int_{\bar y_c(t)}^\infty e^{cy}\Bigl(2v(y,t)v_y(y,t) + c v^2(y,t)
  \Bigr)\d y~.
$$
Given $d > -c$, we have $|2vv_y| \le (c{+}d)v^2 + (c{+}d)^{-1}v_y^2$,
hence
\begin{equation}\label{cdbdd}
  e^{c\bar y_c(t)}\epsilon^2 \,\le\, \int_{\bar y_c(t)}^\infty 
  e^{cy}\Bigl(\frac{1}{c+d}\,v_y^2(y,t) + d v^2(y,t)\Bigr)\d y~.
\end{equation}
If we choose $d$ such that $d(c+d) = \beta_1$ and insert the 
resulting inequality into \reff{lowprelim}, we obtain
\begin{equation}\label{lowbd}
  E_c(t) \,\ge\, e^{c\bar y_c(t)}\Bigl(-\frac{A}{c} + \kappa 
  \epsilon^2\Bigr)~, \quad \hbox{where}\quad 0 \,\le\, 
  \kappa \,\le\, \frac14(c+\sqrt{c^2+4\beta_1})~.
\end{equation}
This estimate shows in particular that the energy $E_c(t)$ is bounded
from below as long as the invasion point $\bar y_c(t)$ is bounded from
above.  Moreover, the lower bound is close to zero if $\bar y_c(t)$ is
large and negative.

\medskip\noindent{\bf Variation of $D_c\,$:}
It follows from \reff{eqv} and \reff{DDcdef} that $D_c(t) = \int_\real 
e^{cy} v_t^2(y,t)\d y$. Differentiating this relation with respect 
to $t$ and integrating by parts, we find
\begin{eqnarray*}
  \frac12\,D_c'(t) &=& \int_{\real} e^{cy}(v_t v_{tt})(y,t)\d y
  \,=\, \int_{\real} e^{cy} v_t (v_{tyy} + c v_{ty} -F''(v)v_t)\d y\\
  &=& -\int_{\real} e^{cy} v_{ty}^2(y,t)\d y \,-\, \int_{\real} e^{cy}
  F''(v(y,t))v_t^2(y,t)\d y~.
\end{eqnarray*}
Take $C_1 > 0$ such that $F''(u) \ge -C_1/2$ whenever $|u| \le B$, 
where $B$ is as in \reff{Bdef}. Then the above relation shows that
\begin{equation}\label{Dcdot}
  D_c'(t) \,\le\, C_1 D_c(t)~, \quad t \ge 0~.
\end{equation}
This differential inequality implies that, if $D_c \in L^1(\real_+)$, 
then $D_c(t) \to 0$ as $t \to \infty$. Since $E_c'(t) = -D_c(t)$
by \reff{ED}, this will be the case as soon as $E_c$ is bounded 
from below. 

\medskip\noindent{\bf Lower bound on $D_c\,$:} Using \reff{DDcdef}
again and integrating by parts, we find
\begin{eqnarray*}\nonumber
  D_c(t) &=& \int_\real e^{cy} (v_{yy} + cv_y - F'(v))^2(y,t) \d y \\
  &=& \int_\real e^{cy} \Bigl(v_{yy}^2 + 2F''(v)v_y^2 + F'(v)^2\Bigr)
  (y,t)\d y~. 
\end{eqnarray*}
We split the integration domain into $(-\infty,\bar Y_c(t))$ and 
$(\bar Y_c(t),+\infty)$. Using \reff{Bdef}, \reff{farfield} and the 
lower bound $F''(v) \ge -C_1/2$, we obtain
$$
  D_c(t) \,\ge\, -\frac{C_1 B^2}{c}\,e^{c\bar Y_c(t)} + 
  \int_{\bar Y_c(t)}^\infty e^{cy}\Bigl(v_{yy}^2(y,t) + 
  2\beta_1 v_y^2(y,t) + \beta_1^2 v^2(y,t)\Bigr) \d y~.
$$
Observe that, for any $y_0 \in \real$, 
\begin{equation}\label{poincare}
  \int_{y_0}^\infty e^{cy} v_y^2\d y \,\ge\, 
  \frac{c^2}4 \int_{y_0}^\infty e^{cy} v^2\d y~, \quad \hbox{and}\quad
  \int_{y_0}^\infty e^{cy} v_{yy}^2\d y \,\ge\, 
  \frac{c^2}4 \int_{y_0}^\infty e^{cy} v_y^2\d y~. 
\end{equation}
Indeed the first inequality is just \reff{cdbdd} with $\epsilon = 0$, 
$d = -c/2$, and $\bar y_c(t)$ replaced by $y_0$, and the second 
inequality is similar. Thus, for any $d \ge 0$ we have
\begin{equation}\label{lowDc}
  D_c(t) \,\ge\, -\frac{C_1 B^2}{c}\,e^{c\bar Y_c(t)} + 
  \int_{\bar Y_c(t)}^\infty e^{cy}\left\{\Bigl(2\beta_1 + \frac{c^2}4
  -d\Bigr) v_y^2 + \Bigl(\beta_1^2 + \frac{dc^2}4\Bigr)
  v^2\right\}(y,t) \d y~.
\end{equation}
In an analogous way we find
\begin{eqnarray}\nonumber
  E_c(t) &=& \int_{-\infty}^{\bar Y_c(t)} e^{cy} 
  \Bigl(\frac12\,v_y^2(y,t) + F(v(y,t))\Bigr) \d y + 
  \int_{\bar Y_c(t)}^\infty e^{cy}\Bigl(\frac12\,v_y^2(y,t) + 
  F(v(y,t))\Bigr) \d y \\ \label{upEc}
  &\le& \frac{K}{c}\,e^{c\bar Y_c(t)} + 
  \int_{\bar Y_c(t)}^\infty e^{cy}\Bigl(\frac12\,v_y^2(y,t) + 
  \frac{\beta_2}2 v^2(y,t)\Bigr) \d y~,
\end{eqnarray}
where $K = (B^2/2) + K'$ and $K' = \sup\{F(u)\,|\,|u|\le B\}$.
If we now combine \reff{lowDc}, \reff{upEc} and choose the 
particular value $d = \beta_2 - (\beta_2-\beta_1)^2/(\beta_2 + 
c^2/4) \ge 0$, we arrive at
\begin{equation}\label{DcEc}
  D_c(t) \,\ge\, \gamma E_c(t) - \frac{C_2}{c}\,e^{c\bar Y_c(t)}~,
  \quad t \ge 0~,
\end{equation}
where
$$
  0 \,<\, \gamma \,\le\, \frac12\,\frac{(c^2+4\beta_1)^2}
  {c^2+4\beta_2}~, \quad \hbox{and} \quad C_2 = C_1 B^2 + \gamma K~.
$$
Inequality \reff{DcEc} means that, if the invasion point $\bar Y_c(t)$ 
is large and negative, the energy dissipation $D_c = -E_c'$ is 
essentially proportional to the energy itself. This gives a 
differential inequality for $E_c(t)$ which, in view of 
Lemma~\ref{Xxcontrol}, can be integrated as follows:
\begin{equation}\label{Ecint}
  E_c(t) \,\le\, e^{-\gamma (t-t_0)}E_c(t_0) + \frac{C_2 T_0}{c}
  \,e^{c(\bar y_c(t_0)+C_0)}~, \quad t \in [t_0,t_0+T_0]~. 
\end{equation}

\begin{remark}\label{constants}
The constants $C_0, C_1$ and $T_0$ introduced in this section depend 
only on the potential $F$. In particular, they are independent of
the solution $u(x,t)$ and of the speed parameter $c$. Similarly, 
if we choose $\kappa = \sqrt{\beta_1}/2$ and $\gamma = 
2\beta_1^2/\beta_2$, then the constants $\kappa, \gamma$ and $C_2$
depend only on $F$. 
\end{remark}

\section{Existence of the invasion speed} \label{invspeed}

As in the previous section, we suppose that $u(x,t)$ is a solution of
\reff{equ} with initial data $u_0$ satisfying the assumptions of
Proposition~\ref{thm2}. We also assume that the bounds \reff{Bdef}
hold for all $t \ge 0$. If $\bar x(t)$ denotes the invasion point
\reff{barxdef}, we define
\begin{equation}\label{c+-def}
c_- \,=\, \liminf_{t\to\infty} \frac{\bar x(t)}{t}~, \qquad
c_+ \,=\, \limsup_{t\to\infty} \frac{\bar x(t)}{t}~.
\end{equation}
Our first result shows that the solution $u(x,t)$ invades 
the stable equilibrium $u = 0$ at a positive, but finite, 
speed. 

\begin{proposition}\label{c+-bounds}
One has $c_- > 0$ and $c_+ < \infty$. 
\end{proposition}

\proof 
The proof relies on the lower bound \reff{lowbd}. Assume that the
initial data $u_0$ belong to $H^1_{c_0}(\real)$ for some $c_0 >
\sqrt{2A}/\epsilon$, where $A = -F(1)$ and $\epsilon$ is as in
\reff{epsdef}. Using \reff{lowbd} with $c = c_0$ and $\kappa = c/2$,
we find that $E_c(t) \ge \alpha \,e^{c\bar y_c(t)}$ for some $\alpha >
0$. Since $E_c(t) \le E_c(0)$ for all $t \ge 0$, it follows that $\bar
y_c(t) = \bar x(t) - ct$ is bounded from above, hence
$$
   c_+ \,=\, c \,+\, \limsup_{t\to\infty} \frac{\bar y_c(t)}{t}
   \,\le\, c \,<\,\infty~.
$$
On the other hand, since $u_0 - 1 \in H^1(\real_-)$ and 
$F(1) = -A < 0$, it is easy to verify that $E_c(0) = \EE_c[u_0]
\sim -A/c$ as $c \to 0$. Thus if we take $c > 0$ sufficiently 
small so that $E_c(0) < 0$, it follows from \reff{lowbd} that 
$0 > E_c(0) \ge E_c(t) \ge (-A/c)\,e^{c\bar y_c(t)}$ for all 
$t \ge 0$. This implies that $\bar y_c(t) = \bar x(t) - ct$ is 
bounded from below, hence 
$$
   c_- \,=\, c \,+\, \liminf_{t\to\infty} \frac{\bar y_c(t)}{t}
   \,\ge\, c \,>\, 0~.
$$
This concludes the proof. \QED

\medskip We next prove that the average invasion speed 
$\bar x(t)/t$ converges to a limit as $t \to \infty$. 

\begin{proposition}\label{cexists}
One has $c_- = c_+$. 
\end{proposition}

\proof We argue by contradiction. Assume that $c_- < c_+$, and choose 
time sequences $\{t_n\}_{n\in\natural}$, $\{t_n'\}_{n\in\natural}$ 
such that $t_n \to \infty$, $t_n' \to \infty$ and 
$$
  \frac{\bar x(t_n')}{t_n'} \,\build\hbox to 8mm
  {\rightarrowfill}_{n \to \infty}^{}\, c_-~, \qquad
  \frac{\bar x(t_n)}{t_n} \,\build\hbox to 8mm
  {\rightarrowfill}_{n \to \infty}^{}\, c_+~.
$$
Due to Lemma~\ref{bas}, upon extracting a subsequence we can assume 
that $u(\bar x(t_n)+z,t_n)$ converges in $H^2_\loc(\real)$ to some 
limit $w_\infty(z)$. More precisely, for any $L > 0$, 
\begin{eqnarray*}
  &&u(\bar x(t_n)+z,t_n) 
  \,\build\hbox to 8mm{\rightarrowfill}_{n \to \infty}^{}\, 
  w_\infty(z) \quad\hbox{in}\quad H^2([-L,L])~, \\
  &&u_t(\bar x(t_n)+z,t_n) 
  \,\build\hbox to 8mm{\rightarrowfill}_{n \to \infty}^{}\, 
  \hat w_\infty(z) \quad\hbox{in}\quad L^2([-L,L])~,
\end{eqnarray*}
where $w_\infty \in H^2_\loc(\real) \cap L^\infty(\real)$ and
$\hat w_\infty \in L^\infty(\real)$ satisfy $\hat w_\infty =
w_\infty'' - F'(w_\infty)$.  Moreover, by definition of the invasion 
point, we have $|w_\infty(0)| = \epsilon$. 

Now, we fix any $c \in (c_-,c_+)$ and we observe that the invasion 
point $\bar y_c(t) = \bar x(t) - ct$ satisfies $\bar y_c(t_n') 
\to -\infty$ and $\bar y_c(t_n) \to +\infty$ as $n \to \infty$. 
Using first the lower bound \reff{lowbd}, we find
\begin{equation}\label{lowbd0}
  E_c(t_n') \,\ge\, -\frac{A}{c}\,e^{c\bar y_c(t_n')} 
  \,\build\hbox to 8mm{\rightarrowfill}_{n \to \infty}^{}\, 0~,
\end{equation}
hence (since $E_c$ is non-increasing) $E_c(t) \ge 0$ for all 
$t \ge 0$. As $E_c'(t) = -D_c(t)$, we deduce that
\begin{equation}\label{Dcint}
  \int_0^\infty D_c(t)\d t \,\le\, E_c(0) \,<\, \infty~,
\end{equation}
and using in addition \reff{Dcdot} we conclude that $D_c(t) \to 0$ 
as $t \to \infty$. 

Next, we observe that, for all $n \in \natural$,
\begin{eqnarray}\nonumber
  D_c(t_n) &=& \int_\real e^{cy} v_t^2(y,t_n)\d y \,=\,
  e^{c\bar y_c(t_n)} \int_\real e^{cz} v_t^2(\bar y_c(t_n) + z,t_n)
  \d z\\ \label{relax}
  &=& e^{c\bar y_c(t_n)} \int_\real e^{cz} (u_t + cu_x)^2
  (\bar x(t_n) + z,t_n)\d z~.
\end{eqnarray}
Since $D_c(t_n) \to 0$ as $n \to \infty$, the last integral in 
\reff{relax} converges to zero as $n \to \infty$, hence the 
limits $w_\infty, \hat w_\infty$ satisfy $\hat w_\infty + 
c w_\infty' = 0$. Incidentally, this means that $w_\infty'' +
cw_\infty' - F'(w_\infty) = 0$, i.e. $w_\infty$ is a travelling wave 
solution of \reff{equ} with speed $c$. Now the crucial point is
that $c \in (c_-,c_+)$ is {\em arbitrary}. Obviously, the relation
$\hat w_\infty + c w_\infty' = 0$ can be satisfied for two different
values of $c$ only if $w_\infty' \equiv 0$, i.e. if $w_\infty$ 
is identically constant. But then we must have $F'(w_\infty) = 0$, 
which is impossible in view of \reff{farfield} since $|w_\infty| = 
\epsilon$. This contradicts the assumption $c_- < c_+$ and 
concludes the proof. \QED

\begin{remark}\label{otherarg}
Another way to obtain a contradiction in the proof of 
Proposition~\ref{cexists}, which works even if $u = 0$ 
is not a strict local minimum of $F$ (see hypothesis {\bf H2}
in the introduction), is to observe that the limiting function 
$w_\infty(z)$ converges to zero as $z \to +\infty$. Indeed, 
proceeding as in \reff{lowprelim}, \reff{relax} and using 
\reff{poincare} we find for all $n \in \natural$:
\begin{eqnarray}\nonumber
  E_c(t_n) &=& e^{c\bar y_c(t_n)} \int_\real e^{cz} 
  \Bigl(\frac12\,u_x^2 + F(u)\Bigr)(\bar x(t_n)+z,t_n)\d z \\
  \label{othereq}  
  &\ge& e^{c\bar y_c(t_n)} \int_0^\infty e^{cz} \Bigl(\frac12\,u_x^2 
  + \frac{\beta_1}2\,u^2\Bigr)(\bar x(t_n)+z,t_n)\d z \,-\, 
  \frac{A}{c}\,e^{c\bar y_c(t_n)} \\ \nonumber
  &\ge& e^{c\bar y_c(t_n)} \,\frac{\kappa'}{2} \int_0^\infty e^{cz} 
  (u_x^2 + u^2)(\bar x(t_n)+z,t_n)\d z \,-\, 
  \frac{A}{c}\,e^{c\bar y_c(t_n)}~,
\end{eqnarray}
where $\kappa' = \min\{1,(c^2{+}4\beta_1)(c^2{+}4)^{-1}\}$. 
As $E_c(t_n) \le E_c(0)$ and $\bar y_c(t_n) \to +\infty$ as 
$n \to \infty$, we have by Fatou's lemma:
\begin{eqnarray*}
  &&\int_0^\infty e^{cz} (w_\infty'(z)^2 + w_\infty(z)^2)\d z\\ 
  &&\qquad \le\, \liminf_{n \to \infty} \int_0^\infty e^{cz} (u_x^2 + 
  u^2)(\bar x(t_n)+z,t_n)\d z \,\le\, \frac{2A}{c\kappa'} 
  \,<\, \infty~. 
\end{eqnarray*}
Thus $w_\infty \in H^1_c(\real)$, and in particular $w_\infty(z) \to
0$ as $z \to +\infty$. This is clearly impossible if $w_\infty' 
\equiv 0$ and $|w_\infty(0)| = \epsilon$. 
\end{remark}

\section{Local convergence to a travelling wave}
\label{localconv}

This section is devoted to the proof of Proposition~\ref{thm2}. 
Using the same notations as in the previous sections, we first
prove that the solution $u(x,t)$ of \reff{equ} converges for 
a {\em sequence of times} towards a travelling wave, locally in
space around the invasion point. On this occasion we identify 
the invasion speed given by Proposition~\ref{cexists} with the 
unique speed $c_*$ for which travelling waves exist. 

\begin{proposition}\label{conv1}
Let $c_\infty = c_- = c_+$. There exists a sequence $t_n \to \infty$
such that, for all $L > 0$, 
\begin{equation}\label{convL}
  \int_{-L}^L e^{c_\infty z}(u_t + c_\infty u_x)^2
  (\bar x(t_n) + z,t_n)\d z 
  \,\build\hbox to 8mm{\rightarrowfill}_{n \to \infty}^{}\, 0~.
\end{equation}
\end{proposition}

\proof Since the left-hand side of \reff{convL} is a nondecreasing
function of $L$, it is sufficient to prove that, for any $L > 0$, 
there exists a sequence $t_n \to \infty$ such that \reff{convL} 
holds. We argue by contradiction and assume that there exist 
$L > 0$ and $\delta > 0$ such that
\begin{equation}\label{notconvL}
  \int_{-L}^L e^{c_\infty z}(u_t + c_\infty u_x)^2(\bar x(t) + z,t)
  \d z \,\ge\, \delta~,
\end{equation}
for all sufficiently large $t$. In fact, upon changing the origin
of time, we can assume that \reff{notconvL} holds for all $t\ge0$.
In analogy with \reff{barydef}, we denote $\bar y(t) = \bar x(t) - 
c_\infty t$. Two situations may occur:

\figurewithtex 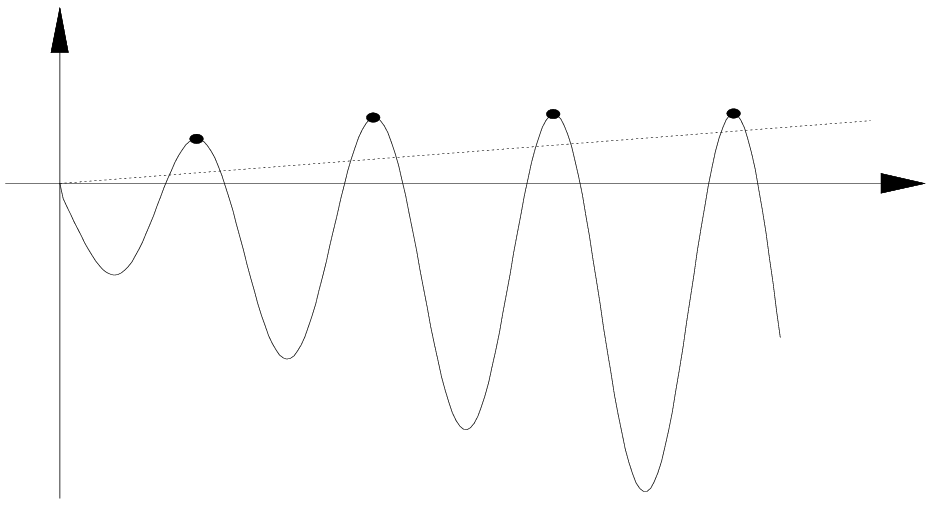 Fig2.tex 5.000 9.000 {\bf Fig.~2:} If there
exists a sequence $t_n \to \infty$ such that $\bar y(t_n)$ is bounded
from below, a contradiction is obtained by considering the dissipation
of the energy $E_c$ in a moving frame with speed $c > c_\infty$ ($c$
close to $c_\infty$). If $\bar y(t_n) \ge 1$ for all $n \in \natural$, 
the set $S_c$ consisting of all $n$ such that $\bar y(t_n) \ge 
(c-c_\infty)t_n$ increases as $c \to c_\infty$, and $\card(S_c) 
\to \infty$.\cr

\medskip\noindent{\bf Case 1:} There exists a time sequence 
$t_n \to \infty$ such that $\{\bar y(t_n)\}_{n\in\natural}$ 
is bounded from below. Without loss of generality we can 
assume that $t_{n+1} \ge t_n + 1$ and $\bar y(t_n) \ge 1$
for all $n \in \natural$ (the second condition is easily achieved
by translating the origin).  

Let $K > 0$ be such that $E_c(0) \le K$ for all $c \in
[c_\infty,c_0]$, where $c_0 > c_\infty$ is as in the proof of
Proposition~\ref{c+-bounds}. Take $c \in (c_\infty,c_0)$ 
sufficiently close to $c_\infty$ so that
\begin{equation}\label{ccond1}
  e^{(c-c_\infty)L} \,\le\, 2~, \quad \hbox{and} \quad 
  (c-c_\infty)^2 B^2 \int_{-L}^L e^{c_\infty z} \d z \,\le\,
  \frac{\delta}{4}~,
\end{equation}
where $B$ is as in \reff{Bdef}. Let $\bar y_c(t) = \bar x(t) - ct = 
\bar y(t) - (c{-}c_\infty)t$. Since $\bar y(t_n) \ge 1$ for all
$n \in \natural$, it is clear that the cardinality of the set
$$
  S_c \,=\, \{n \in \natural\,|\, \bar y_c(t_n) \ge 0\} 
  \,=\, \{n \in \natural\,|\, \bar y(t_n) \ge (c-c_\infty)t_n\} 
$$
becomes arbitrarily large as $c \to c_\infty$, see Fig.~2. 
On the other hand, $\bar y_c(t_n) \to -\infty$ as $n \to \infty$ and 
this implies (as in the proof of Proposition~\ref{cexists}) that 
$E_c(t) \ge 0$ for all $t \ge 0$. But for all $n \in S_c$, we have 
by \reff{relax}, \reff{notconvL}, \reff{ccond1}
\begin{eqnarray*}
  D_c(t_n) &=& e^{c \bar y_c(t_n)} \int_\real e^{cz} (u_t + c u_x)^2
  (\bar x(t_n) + z,t_n)\d z \\ 
  &\ge& \int_{-L}^L e^{c z} (u_t + c u_x)^2(\bar x(t_n) + z,t_n)
  \d z \,\ge\, \frac{\delta}8~,
\end{eqnarray*}
because $e^{cz} \ge \frac12 e^{c_\infty z}$ for $|z| \le L$ and 
$(u_t+cu_x)^2 \ge \frac12(u_t+c_\infty u_x)^2 - (c-c_\infty)^2 B^2$.
Moreover, it follows from \reff{Dcdot} that $D_c(t) \ge D_c(t_n)
\,e^{-C_1}$ for all $t \in [t_n-1,t_n]$, hence
$$
  E_c(t_n{-}1) - E_c(t_n) \,=\, \int_{t_n-1}^{t_n} D_c(t)\d t 
  \,\ge\, \frac{\delta}8\,e^{-C_1}~, \quad n \in S_c~.
$$
If we choose $c$ close enough to $c_\infty$ so that $\card(S_c) > 8K
e^{C_1}/\delta$, we obtain a contradiction with the fact 
that $E_c(t)$ is positive, nonincreasing, and $E_c(0) \le K$. 

\figurewithtex 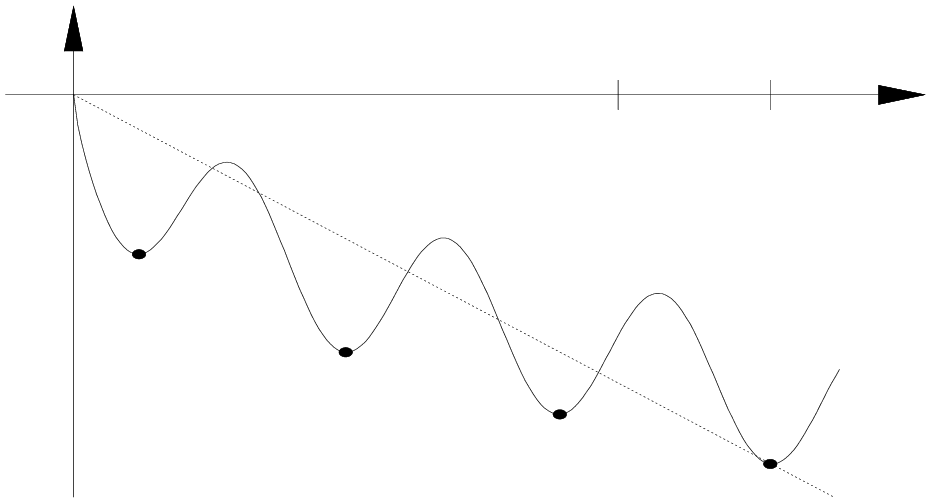 Fig3.tex 5.000 9.000
{\bf Fig.~3:} If $\bar y(t) \to -\infty$ a contradiction is obtained
by considering the dissipation of the energy $E_{c_n}$ in a moving
frame with speed $c_n < c_\infty$ on the time interval $[0,t_n]$,
where $\bar y(t_n) = (c_n-c_\infty)t_n$. We choose $T \gg 1$ and then
$n$ large enough so that $(c_\infty-c_n)T \le 1$.\cr

\medskip\noindent{\bf Case 2:} $\bar y(t) \to -\infty$ as 
$t \to \infty$. In this case, there exists a sequence $t_n \to 
\infty$ such that
\begin{equation}\label{infseq}
  \bar y(t_n) \,\le\, \inf_{0 \le s \le t_n} \bar y(s) + 1~, 
  \quad \hbox{for all } n \in \natural~,
\end{equation}
see Fig.~3. Indeed the function $\mu(t) = \inf\{\bar y(s)\,|\, 
0 \le s \le t\}$ is nonincreasing and $\mu(t) \to -\infty$ as
$t \to \infty$. For each $n \in \natural$, we choose $t_n 
\in [0,n]$ such that $\bar y(t_n) \le \mu(n) + 1$. Then 
$\mu(n) \le \mu(t_n)$, hence \reff{infseq} holds. 

Given some (large) $n \in \natural$, we take $c_n < c_\infty$
such that $(c_n - c_\infty)t_n = \bar y(t_n)$, or equivalently 
$\bar y_{c_n}(t_n) = 0$. Since $c_n \to c_\infty$ as $n \to \infty$,
we can assume that $c_n \ge c_\infty/2$ and 
\begin{equation}\label{ccond2}
  e^{(c_\infty-c_n)L} \,\le\, 2~, \quad 
  (c_\infty-c_n)^2 B^2 \int_{-L}^L e^{c_\infty z} \d z \,\le\,
  \frac{\delta}{4}~.
\end{equation}
If $t \in [0,t_n]$, we have by \reff{infseq}
$$
  \bar y_{c_n}(t) \,=\, \bar y(t) + (c_\infty-c_n)t 
  \,\ge\, \bar y(t_n) - 1 + (c_\infty-c_n)t 
  \,=\, (c_n - c_\infty)(t_n-t)-1~.
$$
Using \reff{notconvL}, \reff{ccond2} and proceeding as in the
previous case, we obtain
\begin{eqnarray*}
  D_{c_n}(t) &=& e^{c_n \bar y_{c_n}(t)} \int_\real e^{c_n z} 
  (u_t + c_n u_x)^2(\bar x(t) + z,t) \d z \\
  &\ge& e^{c_n((c_n-c_\infty)(t_n-t)-1)} \int_{-L}^L e^{c_n z} 
  (u_t + c_n u_x)^2(\bar x(t_n) + z,t_n) \d z\\ 
  &\ge& e^{c_\infty((c_n-c_\infty)(t_n-t)-1)} \,\frac{\delta}8~, 
\end{eqnarray*}
hence for all $T \le t_n$:
\begin{equation}\label{dissipA}
  \int_{t_n-T}^{t_n} D_{c_n}(t)\d t \,\ge\, T \,e^{c_\infty
  ((c_n-c_\infty)T-1)} \,\frac{\delta}8~. 
\end{equation}
On the other hand, there exists $K > 0$ such that $E_{c_n}(0) 
\le K$ for all $n$, and since $\bar y_{c_n}(t_n) = 0$ we know from 
\reff{lowbd} that $E_{c_n}(t_n) \ge -A/c_n$. Thus
\begin{equation}\label{dissipB}
  \int_0^{t_n} D_{c_n}(t)\d t \,=\, E_{c_n}(0) - E_{c_n}(t_n) 
  \,\le\, K + \frac{A}{c_n} \,\le\, K + \frac{2A}{c_\infty}~.
\end{equation}
If we now choose $T > 0$ large enough so that $T e^{-2c_\infty}\delta
> 8K + 16A/c_\infty$, and then $n \in \natural$ large enough so that
$t_n \ge T$ and $(c_\infty-c_n)T \le 1$, we obtain the desired 
contradiction by comparing \reff{dissipA} and \reff{dissipB}.  
\QED

\begin{corollary}\label{conv2}
One has $c_- = c_+ = c_*$, and there exists a sequence $t_n \to \infty$
such that, for all $L > 0$, 
$$
  \sup_{z \in [-L,L]} |u(\bar x(t_n)+z,t_n) - h_\epsilon(z)| 
  \,\build\hbox to 8mm{\rightarrowfill}_{n \to \infty}^{}\, 0~.
$$
\end{corollary}

\proof We argue as in the proof of Proposition~\ref{cexists}.  
If $\{t_n\}_{n \in \natural}$ is the sequence given by
Proposition~\ref{conv1}, we know that (upon extracting a subsequence)
$u(\bar x(t_n)+z,t_n)$ converges in $H^2_\loc(\real)$ to a limit
$w_\infty(z)$ which satisfies
$$
  w_\infty''(z) + c_\infty w_\infty'(z) - F'(w_\infty(z)) \,=\, 0~,
  \quad z \in \real~.
$$
Moreover, $|w_\infty(z)| \le \epsilon$ for all $z \ge 0$,
$|w_\infty(z)| \le B$ for all $z \le 0$, and $|w_\infty(0)| = 
\epsilon$. Arguing as in Remark~\ref{otherarg}, one can also 
show that $w_\infty(z) \to 0$ as $z \to +\infty$. These properties 
together imply that $c_\infty = c_*$ and that $w_\infty = h_\epsilon$, 
see hypothesis {\bf H3} in the introduction. \QED

\begin{corollary}\label{conv3}
For all $c \ge c_*$ and all $w \in H^1_c(\real)$, one has
$$
  \EE_c[w] \,=\, \int_\real e^{cx} \Bigl(\frac12\,w'(x)^2 + 
  F(w(x))\Bigr)\d x \,\ge\, 0~.
$$
\end{corollary}

\proof Assume first that $w \in H^1_{c'}(\real)$ for all $c' > 0$, and
that $w - 1 \in H^1(\real_-)$. If $u(x,t)$ is the solution of
\reff{equ} with initial data $u(x,0) = w(x)$, we know from
Proposition~\ref{cexists} and Corollary~\ref{conv2} that the invasion
point $\bar x(t)$ defined by \reff{barxdef} satisfies $\bar x(t)/t \to
c_*$ as $t \to \infty$.  Thus, for any $c > c_*$, the quantity $\bar
y_c(t) = \bar x(t) -ct$ converges to $-\infty$ as $t \to \infty$, so
that $E_c(t) \ge 0$ for all $t \ge 0$. In particular, $E_c(0) =
\EE_c[w] \ge 0$. Letting $c \to c_*$, we also obtain $\EE_{c_*}[w] \ge
0$.

Assume now that $c \ge c_*$ and that $w \in H^1_c(\real)$. For any 
$n \ge 1$ we define
$$
  w_n(x) \,=\, w(x)\chi(x{-}n) + (1-w(x))\chi(x{+}n{+}1)~, \quad
  x \in \real~,
$$
where $\chi \in \CC^\infty(\real)$, $\chi(x) = 1$ for $x \le 0$ and
$\chi(x) = 0$ for $x \ge 1$. Clearly $w_n(x) = w(x)$ for $x \in
[-n,n]$, whereas $w_n(x) = 0$ for $x \ge n+1$ and $w_n(x) = 1$ for $x
\le -n-1$. Thus $w_n \in H^1_{c'}(\real)$ for all $c' > 0$ and $w_n -
1 \in H^1(\real_-)$, so that $\EE_c[w_n] \ge 0$ for all $n \in
\natural$ by the preceding argument. Moreover it is straightforward to
verify that $\EE_c[w_n] \to \EE_c[w]$ as $n \to \infty$, hence
$\EE_c[w] \ge 0$. \QED

\medskip
Equipped with these results, we are now able to prove that the solution
$u(x,t)$ converges for {\em all times} towards a travelling wave, 
locally in space around the invasion point.

\begin{proposition}\label{conv4}
For all $L > 0$ we have
\begin{equation}\label{convLL}
  \int_{-L}^L e^{c_* z}(u_t + c_* u_x)^2 (\bar x(t) + z,t)\d z 
  \,\build\hbox to 8mm{\rightarrowfill}_{t \to \infty}^{}\, 0~.
\end{equation}
\end{proposition}

\proof We argue by contradiction and assume that there exist $L > 0$, 
$\delta > 0$, and a sequence $t_n \to \infty$ such that
\begin{equation}\label{notconvLL}
  \int_{-L}^L e^{c_* z}(u_t + c_* u_x)^2(\bar x(t_n) + z,t_n)\d z 
  \,\ge\, \delta~,
\end{equation}
for all $n \in \natural$. Let $\bar y(t) = \bar x(t) - c_* t$. 
If the sequence $\{\bar y(t_n)\}_{n \in \natural}$ has a subsequence 
that is bounded from below, then we easily get a contradiction 
as in the proof of Proposition~\ref{conv1} (case 1). 
So it remains to consider the case where $\bar y(t_n) \to -\infty$, 
which requires a new argument. Without loss of generality, we can 
suppose that $t_{n+1} \ge t_n + T_0$ for all $ \in \natural$, where 
$T_0 > 0$ is as in Lemma~\ref{Xxcontrol}, and that $\bar y(t_n) 
\le -n-1$. Upon extracting a subsequence, we can also assume that 
$u(\bar x(t_n)+z,t_n)$ converges in $H^2_\loc(\real)$ towards 
a limit $w_\infty(z)$. 

\figurewithtex 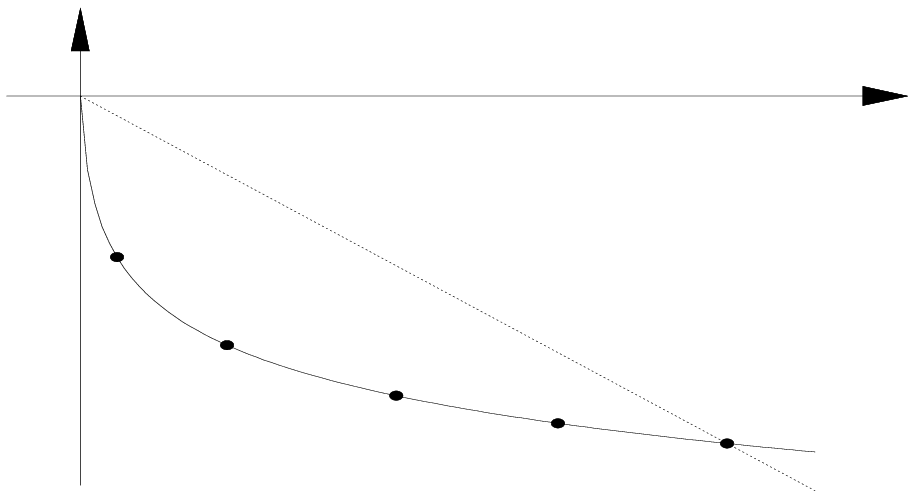 Fig4.tex 5.000 9.000
{\bf Fig.~4:} If $\bar y(t_n) \to -\infty$ a contradiction is obtained
by considering the dissipation of the energy $E_{c_n}$ in a moving
frame with speed $c_n < c_*$ on the time interval $[0,t_n]$,
where $\bar y(t_n) = (c_n-c_*)t_n$.\cr
Given some (large) $n \in \natural$, we take $c_n < c_*$ such 
that $\bar y(t_n) = (c_n-c_*)t_n$, see Fig.~4. Since $c_n \to c_*$
as $n \to \infty$, we can assume that $c_n \ge c_*/2$. Let 
$\bar y_{c_n}(t) = \bar y(t) + (c_*-c_n)t$, so that 
$\bar y_{c_n}(t_n) = 0$. For each $k = 0,1,\dots,n$ we have
by \reff{Ecint}
$$
  E_{c_n}(t_{k+1}) \,\le\, E_{c_n}(t_k + T_0) \,\le\, 
  e^{-\gamma T_0}E_{c_n}(t_k) + 
  \frac{C_2 T_0}{c_n} \,e^{c_n(\bar y_{c_n}(t_k)+C_0)}~, 
$$
hence
\begin{equation}\label{Ecnbdd}
  E_{c_n}(t_k) \,\le\, e^{-k\gamma T_0}E_{c_n}(t_0) + 
  \frac{C_2 T_0 \,e^{c_n C_0}}{c_n} \sum_{j=1}^k e^{-(j-1)\gamma T_0}
  \,e^{c_n\bar y_{c_n}(t_{k-j})}~. 
\end{equation}
We now define $k(n)$ as the largest integer $k \in \natural$ such 
that
$$
  (c_*-c_n)t_j \,\le\, 1 + \frac{j}{2}~, \quad \hbox{for all} 
  \quad j = 0,1,\dots,k~.
$$
Since $c_n \to c_*$, it is clear that $k(n) \to \infty$ as 
$n \to \infty$. Moreover, $k(n) < n$ as $(c_*-c_n)t_n = 
-\bar y(t_n) \ge n+1$ by assumption. For $k = k(n)$ and $j \le k$, 
we have
$$
  \bar y_{c_n}(t_{k-j}) \,=\, \bar y(t_{k-j}) + (c_*-c_n)t_{k-j}
  \,\le\, -(k-j)/2~,
$$
hence it follows from \reff{Ecnbdd} that 
$$
  E_{c_n}(t_n) \,\le\, E_{c_n}(t_{k(n)}) \,\le\, 
  e^{-k(n)\gamma T_0}E_{c_n}(0) +  \frac{C_2 T_0\,e^{c_n C_0}}{c_n} 
  \,k(n)\,e^{-\gamma^*(k(n)-1)}~,
$$
where $\gamma^* = \min(\gamma T_0,c_*/4)$. Taking the limit 
$n \to \infty$ and using the fact that $E_{c_n}(0)$ is 
uniformly bounded, we conclude that
$$
  \limsup_{n \to \infty} E_{c_n}(t_n) \,\le\, 0~.
$$

Now, since $\bar y_{c_n}(t_n) = 0$ by our choice of $c_n$, we 
have
$$
  E_{c_n}(t_n) \,=\, \int_\real e^{c_n z} \Bigl(\frac12 u_x^2
  + F(u)\Bigr)(\bar x(t_n)+z,t_n)\d z~, 
$$
hence taking the limit $n \to \infty$ and using Fatou's lemma
we obtain
$$
  \EE_{c_*}[w_\infty] \,=\, \int_\real e^{c_* z} \Bigl(\frac12 
  w_\infty'(z)^2 + F(w_\infty(z))\Bigr)\d z \,\le\, 
  \liminf_{n \to \infty} E_{c_n}(t_n) \,\le\, 0~.
$$
In particular $w_\infty \in H^1_{c_*}(\real)$, hence it follows 
from Corollary~\ref{conv3} that $\EE_{c_*}[w_\infty] = 0$. 
On the other hand, in view of \reff{notconvLL}, we have
\begin{eqnarray*}
  \DD_{c_*}[w_\infty] &=& \int_\real e^{c_* z} 
  \Bigl(w_\infty''(z) + c_* w_\infty'(z) - F'(w_\infty(z))
  \Bigr)^2\d z\\ 
  &\ge& \int_{-\infty}^L e^{c_* z} \Bigl(w_\infty''(z) + 
  c_* w_\infty'(z) - F'(w_\infty(z))\Bigr)^2\d z\\
  &=& \lim_{n\to\infty} \int_{-\infty}^L e^{c_* z} (u_t + 
  c_* u_x)^2(\bar x(t_n)+z,t_n)\d z \,\ge\, \delta~.
\end{eqnarray*}
Thus, if $u(x,t)$ is the solution of \reff{equ} with initial data
$u(x,0) = w_\infty(x)$, then $E_{c_*}(0) = \EE_{c_*}[w_\infty] = 0$
and $E_{c_*}'(0) = -\DD_{c_*}[w_\infty] \le -\delta$, hence
$E_{c_*}(t) < 0$ for all $t > 0$. This contradicts the conclusion of
Corollary~\ref{conv3}. \QED

\medskip 
It is now a straightforward task to conclude the proof of 
Proposition~\ref{thm2}. Using Proposition~\ref{conv4} and 
proceeding as in Corollary~\ref{conv2}, we see that 
$u(\bar x(t)+z,t)$ converges to $w_\infty(z) \equiv h_\epsilon(z)$ in 
$H^2([-L,L])$ for any $L > 0$. On the other hand, arguing
as in \reff{othereq}, we find for any $c \in (0,c_*)$:
$$
  \limsup_{t \to \infty} \int_0^\infty e^{cz} (u_x^2 + u^2)
  (\bar x(t)+z,t)\d z \,\le\, \frac{2A}{c\kappa'} \,<\, \infty~.
$$
This implies in particular that $u(\bar x(t)+z,t)$ converges 
to zero as $z \to +\infty$ uniformly in $t \ge 0$, hence 
$u(\bar x(t)+z,t)$ converges as $t \to \infty$ to $h_\epsilon(z)$
uniformly for all $z \in [-L,+\infty)$. This proves \reff{locconv}. 

It remains to verify that the map $t \mapsto \bar x(t)$ is $\CC^1$
for large $t$ and satisfies $\bar x'(t) \to c_*$ as $t \to \infty$. 
Using \reff{locconv}, \reff{Bdef}, and an interpolation argument, 
we find for any $L > 0$:
$$
  \sup_{z \in [-L,L]} |u_x(\bar x(t)+z,t) - h_\epsilon'(z)| 
  \,\build\hbox to 8mm{\rightarrowfill}_{t \to \infty}^{}\, 0~.
$$
As $h_\epsilon'(0) < 0$, this implies in particular that 
$u_x(\bar x(t),t)$ is bounded away from zero for $t$ sufficiently
large. Since $u(\bar x(t),t) = \epsilon$ for $t$ large, the 
Implicit Function Theorem then asserts that $\bar x(t)$ is 
differentiable with
$$
  \bar x'(t) \,=\, -\frac{u_t(\bar x(t),t)}{u_x(\bar x(t),t)}~, 
  \quad \hbox{for all sufficiently large } t > 0~.
$$ 
On the other hand $\sup_{|z|\le L} |u_t(\bar x(t){+}z,t) 
+ c_* u_x(\bar x(t){+}z,t)| \to 0$ as $t \to \infty$ by 
\reff{Bdef} and \reff{convLL}, hence $\bar x'(t) \to c_*$ 
as $t \to \infty$. \QED

\section{Repair behind the front}\label{repair}

This final section is devoted to the proof of Corollary~\ref{thm3}.
We follow closely the arguments given in \cite[Section~9.6]{rislerglob},
with a few simplifications. 

Let $u(x,t)$ be a solution of \reff{equ} with initial data $u_0$ 
satisfying the assumptions of Proposition~\ref{thm2}. According to 
\reff{locconv}, we can find a time sequence $t_n \to \infty$ 
such that $t_{n+1} \ge t_n + n + 1$ for all $n \in \natural$, and 
\begin{equation}\label{seqdef}
  \sup_{z \in [-2n,+\infty)} |u(\bar x(t)+z,t) - h_\epsilon(z)|
  \,\le\, \frac{1}{n+1}~, \quad \hbox{for all } t \ge t_n~. 
\end{equation}
Let $\theta : \real \to [0,1]$ be a smooth, nondecreasing function
satisfying $\theta(x) = 0$ for $x \le 0$ and $\theta(x) = 1$ for 
$x \ge 1$. We define a map $\hat x : [0,+\infty) \to \real$ by 
imposing, for all $n \in \natural$:
$$
  \hat x(t) \,=\, \bar x(t) - n - \theta\Bigl(
  \frac{t-t_n}{t_{n+1}-t_n}\Bigr)~, \quad \hbox{for all }
  t \in [t_n,t_{n+1}]~.
$$
It is clear that $\bar x(t)-n-1 \le \hat x(t) \le \bar x(t)-n$ 
for all $t \in [t_n,t_{n+1}]$. Moreover, there exists $T > 0$ 
such that $\hat x(t)$ is differentiable for $t \ge T$, with 
$\hat x'(t) \le \bar x'(t) \le c_* + 1$ for all $t \ge T$. For 
later use we observe that, for any $L > 0$, 
\begin{equation}\label{uuxconv}
  \sup_{z \in [-L,L]}\Bigl(|u(\hat x(t)+z,t)-1| + 
  |u_x(\hat x(t)+z,t)|\Bigr)
  \,\build\hbox to 8mm{\rightarrowfill}_{t \to \infty}^{}\, 0~.
\end{equation}
Indeed, since $\hat x(t) \approx \bar x(t)-n$ for $t \in 
[t_n,t_{n+1}]$, the estimate on $|u-1|$ is a consequence of 
\reff{seqdef} and of the fact that $h_\epsilon(z) \to 1$ as 
$z \to -\infty$. The result for $|u_x|$ then follows from the 
a priori bound \reff{Bdef} by interpolation. 

We next consider the truncated energy function
$$
  E(t) \,=\, \int_\real \phi(x,t)\Bigl(\frac12 u_x^2(x,t) 
  + \overline{F}(u(x,t))\Bigr)\d x~,
$$
where $\overline{F}(u) = F(u) - F(1) \ge 0$ and
$$
  \phi(x,t) \,=\, \left\{
  \begin{array}{ccc}
  1 & \hbox{if} & x \le \hat x(t)~,\\[1mm]
  e^{\hat x(t)-x} & \hbox{if} & x \ge \hat x(t)~.
  \end{array}\right.
$$
Since $u(\cdot,t)-1 \in H^1(\real_-)$ and $u(\cdot,t) \in 
H^1(\real_+)$, it is clear that $E(t)$ is well-defined and 
finite for all $t \ge 0$. Moreover, $E(t)$ is differentiable 
for $t \ge T$ and a direct calculation shows that
\begin{eqnarray*}
  E'(t) &=& -\int_\real \phi(x,t)u_t^2(x,t)\d x 
  + \int_{\hat x(t)}^\infty \phi(x,t)\Bigl\{\hat x'(t) 
  \Bigl(\frac12 u_x^2+\overline{F}(u)\Bigr) + u_x u_t\Bigr\}\d x \\
  &\le& -\frac12 \int_\real \phi(x,t)u_t^2(x,t)\d x  + 
  (c_*{+}1)\int_0^\infty e^{-z}(u_x^2 + \overline{F}(u))
  (\hat x(t)+z,t)\d z~.
\end{eqnarray*}
In view of \reff{Bdef}, \reff{uuxconv}, the last integral in 
the right-hand side converges to zero as $t \to \infty$. 
Since $E(t) \ge 0$ for all $t \ge 0$, it follows that there 
exists a time sequence $t_n' \to \infty$ such that
\begin{equation}\label{tnpdissip}
  \int_\real \phi(x,t_n') u_t^2(x,t_n')\d x 
  \,\build\hbox to 8mm{\rightarrowfill}_{n \to \infty}^{}\, 0~.
\end{equation}

Now, we claim that
\begin{equation}\label{repbehind}
  \sup_{x \in (-\infty,\hat x(t_n')]} |u(x,t_n')-1| 
  \,\build\hbox to 8mm{\rightarrowfill}_{n \to \infty}^{}\, 0~.
\end{equation}
Indeed, if this is not the case, there exist a positive constant
$\epsilon'$, a subsequence $\{t_n''\}_{n\in\natural}$ of
$\{t_n'\}_{n\in\natural}$, and a sequence $\{x_n\}_{n\in\natural}
\subset \real$ such that $x_n \le \hat x(t_n'')$ and $|u(x_n,t_n'')-1|
= \epsilon'$ for all $n \in \natural$.  Without loss of generality, we
can assume that $\epsilon' > 0$ is sufficiently small so that the only
bounded solution of the differential equation $w_{xx} - F'(w) = 0$
with $|w(0)-1| \le \epsilon'$ is $w \equiv 1$, see hypothesis {\bf H4}
in the introduction. In view of \reff{uuxconv}, it is clear that $x_n
- \hat x(t_n'') \to -\infty$ as $n \to \infty$. On the other hand,
upon extracting a subsequence, we can assume that, for all $L > 0$,
\begin{eqnarray*}
  &&u(x_n+z,t_n'') 
  \,\build\hbox to 8mm{\rightarrowfill}_{n \to \infty}^{}\, 
  w_\infty(z) \quad\hbox{in}\quad H^2([-L,L])~, \\
  &&u_t(x_n+z,t_n'') 
  \,\build\hbox to 8mm{\rightarrowfill}_{n \to \infty}^{}\, 
  \hat w_\infty(z) \quad\hbox{in}\quad L^2([-L,L])~,
\end{eqnarray*} 
where $w_\infty \in H^2(\real) \cap L^\infty(\real)$ and 
$\hat w_\infty \in L^\infty(\real)$ satisfy $\hat w_\infty(z) = 
w_\infty''(z) - F'(w_\infty(z))$. However, it follows from 
\reff{tnpdissip} that $\hat w_\infty = 0$, hence $w_\infty $ is a
bounded solution of the differential equation $w_\infty'' - 
F'(w_\infty) = 0$ which satisfies $|w_\infty(0) - 1| = \epsilon'$.
This contradicts the assumption above on $\epsilon'$, hence 
\reff{repbehind} must hold. 

Finally, if we combine \reff{seqdef} and \reff{repbehind}, we 
obtain
$$
  \sup_{z \in \real} |u(\bar x(t_n')+z,t_n')-h_\epsilon(z)| 
  \,\build\hbox to 8mm{\rightarrowfill}_{n \to \infty}^{}\, 0~.
$$
In other words, the solution $u(x,t)$ approaches uniformly on 
$\real$ a translate of the travelling wave $h_\epsilon$ 
for a sequence of times $t_n' \to \infty$. On the other hand, 
the classical results of Sattinger \cite{sattinger} show that, 
if assumptions \reff{min1}, \reff{min0} are satisfied, the 
travelling wave $h$ is {\em asymptotically stable with shift}
in the space $L^\infty(\real)$. In other words, Eq.\reff{conc} holds
for any solution of \reff{equ} which is sufficiently close
(uniformly on $\real$) to a translate of $h$. This is the case 
for $u(\cdot,t_n')$ if $n$ is sufficiently large, hence 
Corollary~\ref{thm3} is proved. \QED

%%%%%%%%%%%%%%%%%%%%%%%%%%%%%%%%%%%%%%%%%%%%%%%%%%%%%%%%%%%%%%%%%%%%%%%%%%%

\bibliographystyle{plain}

\begin{thebibliography}{10}

\bibitem{aronsonweinberger} D. G. Aronson, H. F. Weinberger, 
\textsl{Multidimensional nonlinear diffusion arising in population
genetics}, Adv. Math. {\bf 30} (1978), 33--76.

\bibitem{billinghamneedham} J. Billingham, D.J. Needham:
\textsl{The development of travelling waves in quadra\-tic and cubic
autocatalysis with unequal diffusion rates. I. Permanent form 
travelling waves}, Phil. Trans. R. Soc. Lond. A {\bf 334} (1991), 
1--24. 

\bibitem{fife} P. Fife, 
\textsl{Long time behavior of solutions of bistable nonlinear
diffusion equations}, Arch. Rat. Mech. Anal. {\bf 70} (1979), 
31--46.

\bibitem{fifemcleod1} P. Fife, J. B. McLeod, 
\textsl{The approach of solutions of nonlinear diffusion equations 
to travelling front solutions}, Arch. Rat. Mech. Anal. {\bf 65}
(1977), 335--361.

\bibitem{fifemcleod2} P. Fife, J. B. McLeod, 
\textsl{A phase plane discussion of convergence to travelling fronts
for nonlinear diffusion}, Arch. Rat. Mech. Anal. {\bf 75} (1981), 
281--314.

\bibitem{fisher} R.A. Fisher:
\textsl{The Advance of Advantageous Genes}, Ann. of 
Eugenics {\bf 7} (1937), 355--369. 

\bibitem{gallayjoly} Th. Gallay and R. Joly, in preparation. 

\bibitem{henry} D. Henry, 
\textsl{Geometric theory of semilinear parabolic equations},
Springer-Verlag, Berlin, 1981.

\bibitem{kanel1} Ja.I. Kanel':
\textsl{Stabilization of solutions of the Cauchy problem for equations 
encountered in combustion theory}, Mat. Sbornik (N.S.) {\bf 59} 
(1962), 245--288. 

\bibitem{kanel2} Ja.I. Kanel':
\textsl{Stabilization of the solutions of the equations of combustion 
theory with finite initial functions}, Mat. Sbornik {\bf 65} (1964), 
398--413.

\bibitem{kpp} A.N. Kolmogorov, I.G. Petrovskii, N.S. Piskunov, 
\textsl{Etude de la diffusion avec croissance de la quantit\'e de 
mati\`ere et son application \`a un probl\`eme biologique}, 
Moscow Univ. Math. Bull. {\bf 1} (1937), 1--25. 

\bibitem{muratov} C.B. Muratov,
\textsl{A global variational structure and propagation of disturbances 
in reaction-diffusion systems of gradient type}, 
Discrete Contin. Dyn. Syst. Ser. B {\bf 4} (2004), 867--892. 

\bibitem{rislerglob} E. Risler,
\textsl{Global convergence towards travelling fronts in nonlinear
parabolic systems with a gradient structure}, to appear in 
Ann. Inst. H. Poincar\'e.

\bibitem{roquejoffre1} J.-M. Roquejoffre, 
\textsl{Convergence to travelling waves for solutions of a class of 
semilinear parabolic equations}, J. Differential Equations 
{\bf 108} (1994), 262--295.

\bibitem{RTV} J.-M. Roquejoffre, D. Terman and V. Volpert,
\textsl{Global stability of traveling fronts and convergence towards 
stacked families of waves in monotone parabolic systems}, SIAM J. 
Math. Anal. {\bf 27} (1996), 1261--1269.

\bibitem{roquejoffre2} J.-M. Roquejoffre, 
\textsl{Eventual monotonicity and convergence to travelling fronts 
for the solutions of parabolic equations in cylinders}, Ann. Inst. 
H. Poincar\'e Anal. Non Lin\'eaire {\bf 14} (1997), 499--552. 

\bibitem{sattinger} D. H. Sattinger,
\textsl{On the Stability of Waves of Nonlinear Parabolic Systems}, 
Adv. Math. {\bf 22} (1976), 312--355. 

\bibitem{volpert} A. I. Volpert, V. A. Volpert, V. A. Volpert,
\textsl{Traveling wave solutions of parabolic systems}, Translations 
of Mathematical Monographs {\bf 140}, AMS Providence, 1994. 

\end{thebibliography}

\end{document}